\documentclass{amsart}
\usepackage{amssymb}
\usepackage{amsfonts}
\usepackage{latexsym}

\newtheorem{theorem}{Theorem}[section]
\newtheorem{lemma}[theorem]{Lemma}
\newtheorem{proposition}[theorem]{Proposition}
\newtheorem{corollary}[theorem]{Corollary}
\theoremstyle{definition}
\newtheorem{definition}[theorem]{Definition}
\newtheorem{example}[theorem]{Example}

\newtheorem{notation}[theorem]{Notation}
\newtheorem{remark}[theorem]{Remark}

\newcommand{\C}{{\mathbb C}}

\newcommand{\id}{\mbox{id}}

\newcommand{\End}{\mbox{End}}

\newcommand{\Hom}{\mbox{Hom}}
\newcommand{\Ad}{\mbox{Ad}}

\newcommand{\Ind}{\mbox{Ind}}
\newcommand{\Rep}{\mbox{Rep}}

\newcommand{\Prim}{\mbox{Prim}}
\newcommand{\Perm}{\mbox{Perm}}

\newcommand{\eps}{\varepsilon}

\newcommand{\R}{{\mathcal R}}
\newcommand{\bR}{{\bar \R}}

\newcommand{\bJ}{{\bar J}}

\newcommand{\g}{\mathfrak{g}}

\newcommand{\n}{\mathfrak{n}}

\renewcommand{\a}{\mathfrak{a}}

\begin{document}

\title{Vertex-IRF transformations and quantization of 
dynamical $r$-matrices }

\author{Pavel Etingof}
\address{Etingof : Columbia University,
Department of Mathematics,
2990 Broadway, New York, NY 10027 and
MIT, Department of Mathematics,
77 Massachusetts Avenue, Cambridge, MA 02139-4307}

\author{Dmitri Nikshych}
\address{Nikshych : UCLA, Department of Mathematics, 405 Hilgard Avenue,
Los Angeles, CA 90095-1555 and MIT, Department of Mathematics,
77 Massachusetts Avenue, Room 2-130, Cambridge, MA 02139-4307}

\date{March 13, 2001}
\begin{abstract}
Motivated by the correspondence between the vertex and IRF models
in statistical mechanics,
we define and study a notion of vertex-IRF transformation for
dynamical twists that generalizes a usual gauge transformation.
We use vertex-IRF transformations to quantize all completely
degenerate  dynamical $r$-matrices on finite-dimensional Lie
algebras.
\end{abstract}
\maketitle


\begin{section}
{Introduction}

This note has two goals. One is to describe the notion of a  
vertex-IRF transformation, which appeared in physics (as a map
between the vertex and IRF models, see e.g., \cite{Has})
and has so far remained
a part of folklore, but which, we feel, is worthy of a systematic study.
The other is to give examples of such transformations
leading to quantization of dynamical $r$-matrices which have not been
quantized previously.

Recall that the study of the dynamical quantum Yang-Baxter equation of
Felder \cite{F} leads one to a useful notion of a dynamical twist
in a Hopf algebra (see \cite{BBB, EV, ES2}). 
Namely, such a twist in a quasitriangular
Hopf algebra gives rise to a solution of Felder's equation.
There is a notion of a gauge transformation and gauge equivalence
of dynamical twists; and twists which are gauge equivalent can be, in a sense,
regarded as ``the same'' (e.g., the weak Hopf algebras associated to such
twists as in [EN1] are isomorphic).

A vertex-IRF transformation is a generalization
of a gauge transformation. More precisely, it is structurally similar
to an ordinary gauge transformation, but unlike it, allows
one to turn a known dynamical twist into an essentially new one.
For example, one  can apply vertex-IRF transformations to the identity
twist $1\otimes 1$ and obtain new dynamical twists. It turns out that in 
this way one can quantize a large class 
of dynamical $r$-matrices, which we call ``completely degenerate''.
In some sense they are opposite to the non-degenerate dynamical 
$r$-matrices, which were quantized in \cite{X}.  

The structure of the paper is as follows.
 
In Section 2 we define vertex-IRF transformations 
(with discrete and continuous dynamical variable), 
study their main properties, and then define 
vertex-IRF functors, which are a generalization of
the functor between representation categories 
of the Felder and Belavin elliptic quantum groups 
discussed in \cite{Has},\cite{ES1}. 

In Section 3 we give examples of vertex-IRF transforms
and resulting dynamical twists.

In Section 4 we introduce the notion of a completely 
degenerate dynamical $r$-matrix, 
and show that the methods of Section 3 allow one 
to quantize any such $r$-matrix. These quantizations 
appear to be new even in the simplest cases. 

\textbf{Acknowledgments.} We are grateful to Eric Buffenoir,
Phillippe Roche, and Ping Xu for useful discussions. 
In particular, we are grateful to  Ping Xu for asking the question that
led us to the results of Sections 3 and 4. The second author
thanks MIT for the warm hospitality during his visit.
The authors were partially supported by the NSF grant 
DMS-9988796.

\end{section}


\begin{section}
{Vertex-IRF transformations}

\subsection{Definitions of vertex-IRF and IRF-vertex transformations}

Let us give a rigorous definition of a   vertex-IRF
transformation. We will first give the definition for
finite groups, and then generalize it to Lie groups.

We recall definitions of dynamical twists in Hopf algebras
and their gauge equivalence from \cite{EN2}. Let $H$ be a Hopf algebra over 
$\mathbb{C}$ and $A \subset G(H)$ be a finite Abelian subgroup 
of the group $G(H)$ of group-like elements of $H$.

\begin{notation}
\label{h}
Let $P_\mu,\,\mu\in A^*,$ be the set of primitive idempotents in $\C[A]$.
Let $F : A^*\to H^{\otimes n}$ be a function. Set
\begin{equation}
F(\lambda \pm h^{(i)}) = \sum_\mu\, F(\lambda \pm \mu) P_\mu^i
\quad \mbox{ and } \quad
F(\lambda \pm \hat{h}^{(i)}) = \sum_\mu\, P_\mu^i F(\lambda \pm \mu) ,
\end{equation}
where $P_\mu^i = 1\otimes\cdots\otimes P_\mu \otimes\cdots\otimes 1$
with $P_\mu$ in the $i$th component.
\end{notation}

Of course, $F(\lambda \pm h^{(i)}) = F(\lambda \pm \hat{h}^{(i)})$
if $F$ is of zero weight with respect to the action of $A$
in the $i$th component.

\begin{definition}
\label{dynamical twist}
Let $J : A^* \to H\otimes H$ be a
zero weight function with invertible values. We say that  $J$ is a 
{\em dynamical twist}  in $H$ if it satisfies the following 
functional equations 
\begin{eqnarray}
\label{dynamical equation}
&J^{12,3}(\lambda) J^{12}(\lambda-h^{(3)})
= J^{1,23}(\lambda) J^{23}(\lambda),& \\
\label{counit of J}
&(\eps\otimes \id)J(\lambda) = (\id \otimes \eps)J(\lambda) = 1.&
\end{eqnarray}
Here $J^{12,3}(\lambda) =(\Delta\otimes\id)J(\lambda)$,
$J^{12}(\lambda) = J(\lambda)\otimes 1$ etc.
\end{definition}

If $J$ is a dynamical twist in $H$ and $x : A^*\to H$ is a
zero weight function with invertible values
such that $\eps(x(\lambda))\equiv 1$, then 
\begin{equation}
J^x(\lambda) = \Delta(x(\lambda))\,J(\lambda)\, 
                x^1(\lambda-h^{(2)})^{-1} x^2(\lambda)^{-1}
\end{equation}
is also a dynamical twist in $H$.

\begin{definition}
\label{dynamical gauge}
We say that $J^x$ is  {\em gauge equivalent} to $J$ and 
that $x$ is a {\em gauge transformation}.
\end{definition}

Let $\bar A$ be a subgroup of $A$,
$\bar J: \bar A^*\to H\otimes H$ be a dynamical twist, and
$x: A^*\to H$ be a function with invertible values,
of zero weight with respect to $\bar A$ (but not necessarily with 
respect to $A$), such that $\eps(x(\lambda)) \equiv 1$.
Define a function $J^x: A^*\to H\otimes H$ by the formula
\begin{equation}
\label{from bar A to A}
J^x(\lambda)=\Delta(x(\lambda))\,\bar J(\bar\lambda)\,
x^2(\lambda)^{-1} \, x^1(\lambda-h^{(2)})^{-1},
\end{equation}
where $\lambda \mapsto \bar\lambda$ 
is the canonical projection $A^* \to {\bar A}^*$.
Note that the order in which $x^1$ and $x^2$ are written is essential.
Of course, $J^x(\lambda)$ is of zero weight with respect to $\bar A$,
but a priori not  with respect to $A$.

\begin{proposition} 
\label{vertex-IRF}
If $J^x(\lambda)$ is of zero weight with respect to $A$,
then it is a dynamical twist.
\end{proposition}
\begin{proof}
For all $\lambda\in A^*$ we compute, using the zero weight properties
of $x(\lambda)$ and  $J^x(\lambda)$ and the 
dynamical twist equation for  $\bJ(\bar\lambda)$ :
\begin{eqnarray*}
\lefteqn{ {J^x}^{12,3}(\lambda)\, {J^x}^{12}(\lambda-h^{(3)}) = }\\
&=& x^{123}(\lambda)\, \bJ^{12,3}(\bar\lambda)\, x^3(\lambda)^{-1}\,
    \bJ^{12}(\bar\lambda -\bar{h}^{(3)})\, x^2(\lambda-h^{(3)})^{-1}\,
     x^1(\lambda-h^{(2)}-h^{(3)})^{-1} \\
&=& x^{123}(\lambda) \, \bJ^{12,3}(\bar\lambda) \,
    \bJ^{12}(\bar\lambda -\bar{h}^{(3)})\,
    x^3(\lambda)^{-1} \,x^2(\lambda-h^{(3)})^{-1}\, 
    x^1(\lambda-h^{(2)}-h^{(3)})^{-1},\\
\lefteqn{ {J^x}^{1,23}(\lambda)\, {J^x}^{23}(\lambda) = }\\
&=& x^{123}(\lambda)\,  \bJ^{1,23}(\bar\lambda)  \,
    x^{23}(\lambda)^{-1}\, x^1(\lambda-h^{(2)}-h^{(3)})^{-1}\, 
    {J^x}^{23}(\lambda)\\
&=&  x^{123}(\lambda) \, \bJ^{1,23}(\bar\lambda) \, \bJ^{23}(\bar\lambda) \,
     x^3(\lambda)^{-1}\, x^2(\lambda-h^{(3)})^{-1} \,
     x^1(\lambda-h^{(2)}-h^{(3)})^{-1}. 
\end{eqnarray*}
The verification of  equation (\ref{counit of J}) is straightforward.
\end{proof}

\begin{definition} 
\label{vertex-IRF transformation}
If the assumption of Proposition \ref{vertex-IRF}
is satisfied, then the function $x$ is called a vertex-IRF transformation
from $\bar J$ to $J^x$.
\end{definition}

Proposition~\ref{vertex-IRF} has the following ``converse'' version.

Let $J: A^*\to H\otimes H$ be a dynamical twist, and
$x: A^*\to H$ be a function with invertible values, of zero
weight with respect to $\bar A$ (but not necessarily with respect to $A$), 
such that $\eps(x)=1$.
Let the function $J^x: A^*\to H\otimes H$ be given by the formula
\begin{equation}
\label{from A to bar A}
J^x(\lambda)=\Delta(x(\lambda))\,J(\lambda)\,
x^1(\lambda-h^{(2)})^{-1}\,x^2(\lambda)^{-1}.
\end{equation}

\begin{proposition} 
\label{IRF-vertex}
If for some function $\bar J^x : \bar A^*\to H \otimes H$
we have $J^x(\lambda)=\bar J^x(\bar\lambda)$ for all $\lambda\in A^*$,
then $\bar J^x$ is a dynamical twist.
\end{proposition}
\begin{proof}
The computation below is similar to that of Proposition~\ref{vertex-IRF}
and uses the zero weight properties of the functions involved :
\begin{eqnarray*}
\lefteqn{\bJ^{x\,12,3}(\bar\lambda)\,\bJ^{x\,12}(\bar\lambda-\bar{h}^{(3)})=}\\
&=& x^{123}(\lambda)\, {J^x}^{12,3}(\lambda)\, x^{12}(\lambda - h^{(3)})^{-1}\,
     \bJ^{x\,12}(\bar\lambda -\bar{h}^{(3)})\,  x^3(\lambda)^{-1}\\
&=& x^{123}(\lambda)\, {J^x}^{12,3}(\lambda)\, {J^x}^{12}(\lambda -h^{(3)})\,
     x^1(\lambda-h^{(2)}-h^{(3)})^{-1}\, x^2(\lambda -h^{(3)})^{-1} \,
     x^3(\lambda)^{-1}, \\
\lefteqn{\bJ^{x\,1,23}(\lambda)\, \bJ^{x\,23}(\lambda) = }\\
&=& x^{123}(\lambda)\,  {J^x}^{1,23}(\lambda) \,
    x^1(\lambda-h^{(2)}-h^{(3)})^{-1}\, {J^x}^{23}(\lambda)\,
    x^2(\lambda -h^{(3)})^{-1}\, x^3(\lambda)^{-1} \\
&=& x^{123}(\lambda) \, {J^x}^{1,23}(\lambda) \,{J^x}^{23}(\lambda)\,
    x^1(\lambda-h^{(2)}-h^{(3)})^{-1}  \,
    x^2(\lambda -h^{(3)})^{-1} \,x^3(\lambda)^{-1}. 
\end{eqnarray*}
The counit identity is straightforward.
\end{proof}

\begin{definition} 
\label{IRF-vertex transformation}
If the assumption of Proposition \ref{IRF-vertex}
is satisfied, then $x$ is called an   IRF-vertex transformation
from $J$ to $\bJ^x$.
\end{definition}

\begin{remark}
\begin{enumerate}
\item[(1)]
Let $J,\bJ$ be dynamical twists defined on $A^*,\bar A^*$,
respectively. Then $x: A^*\to H$ is a  vertex-IRF transformation from
$\bJ$ to $J$ if and only if $x^{-1}$ is an IRF-vertex transformation from
$J$ to $\bJ$.
\item[(2)]
If $A=\bar A$ then the assumptions of both propositions are vacuous,
and a vertex-IRF or an IRF-vertex transformation is just a usual gauge 
transformation.
\end{enumerate}
\end{remark}

A useful special case of Proposition~\ref{vertex-IRF} is
$\bar{A} =\{ 0\}$ when the twist to be transformed is
$1\otimes 1$. In this case, we have

\begin{corollary}
\label{baratriv} 
Let $x: A^*\to H$ be any invertible-valued function.
If 
\begin{equation*}
J(\lambda):=\Delta(x(\lambda)) x^2(\lambda)^{-1}x^1(\lambda-h^{(2)})^{-1}
\end{equation*}
is of zero weight then it is a dynamical twist.
\end{corollary}

Let us also describe how   vertex-IRF and IRF-vertex transformations
act on $R$-matrices. For this purpose, assume that $H$ is a quasitriangular
Hopf algebra with the universal $R$-matrix ${\mathcal R}$.
In this case, every dynamical twist $J(\lambda)$ defines
the dynamical $R$-matrix 
$\R(\lambda)=J^{21}(\lambda)^{-1}{\mathcal R}J(\lambda)$,
which satisfies the quantum dynamical Yang-Baxter equation.

Let now $J,\,\bJ$ be the twists on $A,\bar A$,
and $x$ a   vertex-IRF transformation such that $\bJ^x=J$.
Let $\R,\,\bR$ be the dynamical $R$-matrices associated
to $J,\,\bJ$.

\begin{corollary}
\label{rmat}
The $R$-matrices above are related by
\begin{equation}
\label{R matrices}
\R(\lambda)=x^2(\lambda-h^{(1)})x^1(\lambda)\,
            \bR(\bar \lambda) \, x^{-2}(\lambda)x^{-1}(\lambda-h^{(2)}).
\end{equation}
\end{corollary}

Let us now extend the above theory to the case of Hopf algebras
over ${\mathbb C}[[\hbar]]$. Let
$H$ be a deformation Hopf algebra over ${\C}[[\hbar]]$ (for example,
a quantized universal enveloping algebra), 
and let $H_0=H/\hbar H$. Let ${\mathfrak a}$ be a
finite-dimensional (over $\mathbb{C}$) commutative Lie subalgebra of
$\Prim(H)$, the Lie algebra of primitive elements of $H$,
such that the induced map ${\mathfrak a}\to \Prim(H_0)$ is
injective. In this situation, we make the following definition.

\begin{definition}
\label{formal dynamical twist} 
A zero weight  $H\otimes H$-valued
meromorphic function $J(\lambda)$ on ${{\mathfrak a}^*}$ is called a {\em
formal dynamical twist } for $H$ if it equals $1\otimes 1$ modulo $\hbar$
and satisfies the following functional equations :
\begin{eqnarray}
\label{formal dynamical equation} 
&J^{12,3}(\lambda) J^{12}(\lambda-\hbar h^{(3)})
= J^{1,23}(\lambda) J^{23}(\lambda),& \\
\label{formal counit of J} 
&(\eps\otimes \id)J(\lambda) = (\id \otimes \eps)J(\lambda)  = 1.&
\end{eqnarray}
\end{definition}

Here the expression $J^{12}(\lambda-\hbar h^{(3)})$ is
understood in the sense of the Taylor expansion with respect to
$\hbar$. 

Note that this definition differs from Definition~\ref{dynamical twist}
only by replacing $h^{(3)}$ with $\hbar h^{(3)}$. Because of this, 
in the sequel we will drop the word ``formal'', and refer to $J$ as
a ``dynamical twist''.

To define vertex-IRF and IRF-vertex transformations
in the formal situation, 
we will define the group $\tilde H^\times$, which is an
extension by the reduced multiplicative group $1+\hbar H$ of $H$
of the simply connected Lie group $G$ 
corresponding to the Lie algebra $\g$ of primitive elements of $H_0$.
To make this definition, we will assume that 

\begin{enumerate}
\item[(1)] the Lie algebra $\g$ is finite-dimensional, and
\item[(2)] the $\g$ module $H_0$ 
(where the action is by the commutator) is a sum of finite-dimensional
submodules. 
\end{enumerate}

These conditions are satisfied, for example, if $H_0$ is the
enveloping algebra of a finite-dimensional Lie algebra. 

Consider the Lie algebra ${\mathfrak l}$ of all elements 
$y\in H$ such that $y\in \g \mod \hbar$
(under commutator). 
We have an exact sequence of Lie algebras 
$$
0\to \hbar H\to {\mathfrak l}\to \g\to 0.
$$
It is easy to show that under assumptions (1) and (2), the Lie
algebra
${\mathfrak l}/\hbar^n {\mathfrak l}$ is a sum of finite-dimensional Lie
algebras. This implies that the above exact sequence of Lie
algebras canonically defines an exact sequence of groups 
$$
0\to 1+\hbar H[[\hbar]]\to L\to G\to 0.
$$ 

This defines the desired group $L$, which we will denote by 
$\tilde H^\times$. We have an exponential map $\exp: {\mathfrak
  l}\to L$, which allows us to think of ${\mathfrak l}$ as a Lie
algebra of $L$. In other words, one may think  
of (small) elements of $\tilde H^\times$ as expressions
of the form $e^F$, where $F\in H$ and $F\mod \hbar \in \g$ 
(in these terms, the multiplication in $\tilde H^\times$ uses the 
Campbell-Hausdorff formula).

Now let $\bar{\mathfrak a}$ be a Lie subalgebra of ${\mathfrak a}$.

We will make the following modifications in the above theory
of vertex-IRF transformations:

\begin{enumerate}
\item[1.] 
$A,\bar A$ are replaced by ${\mathfrak a},\bar {\mathfrak a}$.
\item[2.] 
The functions $x(\lambda)$ are required to be meromorphic and 
take values in $\tilde H^\times$.
\item[3.] 
In all formulas $h^{(i)}$ is replaced with $\hbar h^{(i)}$.
\end{enumerate}

Then we have

\begin{proposition}
Propositions \ref{vertex-IRF}, \ref{IRF-vertex} and Corollaries
\ref{baratriv}, \ref{rmat} are valid in the formal case with the
above modifications.
\end{proposition}

\begin{remark}
\label{Now comes the time}
Let us  motivate the terminology ``vertex-IRF'' and
``IRF-vertex'' transformation. The point is that formula (\ref{R matrices})
is (up to small modifications) exactly the same as the formula
representing the relationship between the $R$-matrices of the 8-vertex
and the interaction-round-a-face models of statistical mechanics,
discovered by Baxter (for $n=2$) and developed by
Hasegawa \cite{Has} (see also \cite[Lemma 2]{ES1}).
In fact, it is known by now (see \cite{JOKS}) that these $R$-matrices
(namely, the Baxter-Belavin and the Felder $R$-matrices in the vector
representation) can be obtained by specializing universal $R$-matrices obtained
from the usual trigonometric $R$-matrix of the quantum affine 
algebra $U_q(\widehat{sl(n)})$
by twisting using a non-dynamical and a dynamical twist, respectively.
It is expected that these two twists  in fact are related by a
 vertex-IRF transformation in $U_q(\widehat{sl(n)})$;
however, the universal expression for this transformation is not known:
it is only known in the vector representation,
and in a tensor product of shifted vector representations (see \cite{ES1}).
The universal expression is known, however, for $n=2$ in a trigonometric
degeneration (see \cite{BBB}, the element $M(\lambda)$). The case $n>2$
was intensively studied by E.~Buffenoir and Ph.~Roche (unpublished),
who turned the authors' attention to the notion of the   vertex-IRF
transformation.  
We also note that the vertex-IRF transformation  for $gl(n)$ 
(in the vector  representation), was, in effect, used by Cremmer and 
Gervais in their proof that the Cremmer-Gervais $R$-matrix satisfies the 
quantum Yang-Baxter equation. This was explained by Hodges 
in \cite{Ho}. 
\end{remark}

\subsection{An IRF-vertex functor}

In this subsection we would like to explain how a vertex-IRF 
transformation  between two dynamical twists gives rise to a certain 
functor between the corresponding categories of representations. 
For the original vertex-IRF transformation coming from physics, 
this functor was essentially introduced by Hasegawa (see \cite{Has}) 
and was studied in \cite{ES1}. 

For simplicity
we will deal with the situation when the dynamical twist is defined on the 
dual of a finite abelian group $A$ lying inside a finite-dimensional Hopf 
algebra $H$; 
the cases of a finite-dimensional Lie algebra
$\a$ instead of $A$ and of infinite dimensional $H$ are completely parallel.

Let $J:A^*\to H\otimes H$ be a dynamical twist. Then one can define
the category $\Rep(J)$ of representations of the dynamical quantum 
group associated to $J$. This can be done, for instance, using weak 
Hopf algebras introduced in \cite{BNSz}.
Namely, one can define a weak Hopf algebra $H^J$ 
($H$ twisted by $J$) as in \cite{EN1}, and 
$\Rep(J)$ is nothing but the category of its comodules.

Let us also  
give another, more explicit description of $\Rep(J)$, 
which can be found in \cite{EV,ES1}. Apart from 
being more explicit, this description has the advantage that, 
unlike the weak Hopf algebra description, it can be easily adapted
to the case when $A$ is replaced with an abelian Lie algebra $\a$. 

An object of $\Rep(J)$ is an
$A$-module $V$, together with an assignment $X\to L_X(\lambda)$, 
which assigns to any $H$-module $X$ 
the L-operator of this module, $L_X:A^*\to \End_A(X\otimes V)$. 
This assignment is required to be compatible with morphisms 
in $\Rep(H)$ and to be compatible with tensor products in the following 
way: 
\begin{equation}
J^{12}(\lambda-h^{(3)})L_Y^{23}(\lambda)
L_X^{13}(\lambda-h^{(2)})J^{12}(\lambda)^{-1} = L_{X\otimes Y}.
\end{equation}
in $X\otimes Y\otimes V$.

A morphism between representations $(V,L),(V',L')$ is a function
$f:A^*\to \Hom_A(V,V')$ such that 
$$
(1\otimes f(\lambda))L_X(\lambda)=L_X'(\lambda)(1\otimes f(\lambda-h^{(1)})).
$$

Now we proceed to construct the desired functor.
Let $A,J$ be as above and  $\bar J\in H\otimes H$ be an ordinary twist.
Assume that $x:A^*\to H$ is a vertex-IRF transformation from $\bar J$
to $J$. 

Let $(V,L)\in \Rep(J)$. For any $X\in \Rep(H)$ define 
the linear operator $\bar L_X$ 
on the vector space $X\otimes 
V\otimes F(A^*)$, where $F(A^*)$ denotes the function 
algebra on $A^*$ : 
$$
\bar L_X:=x^1(\lambda-h^{(2)})^{-1} \circ L_X(\lambda) \circ
T_1  \circ x^1(\lambda+\hat{h}^{(1)}).
$$
Here we use Notation~\ref{h} and set
$T_1 f(\lambda)x\otimes v := f(\lambda-h^{(1)})x\otimes v$.

The main result of this subsection is 

\begin{theorem}
\label{our main result}
\begin{enumerate}
\item[(i)] 
The operators $\bar L_X$ turn the vector space $V\otimes F(A^*)$ into an 
object of  $\Rep(\bar J)$ (i.e., a comodule over the Hopf algebra 
$H^{\bar J}$).
\item[(ii)] 
The assignment ${\mathcal F}:  (V,L)\mapsto (V\otimes F(A^*),
\bar L)$ is a functor.
\end{enumerate}
\end{theorem}

\begin{remark}
We note that here we have considered only the situation $\bar A=0$. 
If $\bar A\ne 0$ and $A=\bar A\oplus \bar A'$ then 
one can generalize the theorem and construct a functor 
$\mathcal{F}$ which, at the level of vector spaces, 
 reduces to tensoring with  $F((\bar{A}')^*)$. 
\end{remark}
\begin{proof} The statements are proved by direct verification.
We will only prove the more difficult statement (i). 
One has
\begin{eqnarray*}
\bar L_{X\otimes Y} 
&=& x^{12}(\lambda-h^{(3)})^{-1}L_{X\otimes Y}(\lambda) 
    T_1 T_2 x^{12}(\lambda+\hat h^{(1)}+\hat h^{(2)})\\
&=& x^{12}(\lambda-h^{(3)})^{-1}J^{12}(\lambda-h^{(3)})L_Y^{23}(\lambda)
L_X^{13}(\lambda-h^{(2)})J^{12}(\lambda)^{-1}  \\
& & \qquad \times\,  T_1 T_2 x^{12}(\lambda+\hat h^{(1)}+\hat h^{(2)})\\
&=& \bar J^{12}x^{2}(\lambda-h^{(3)})^{-1}x^1(\lambda-h^{(2)}-h^{(3)})^{-1}
    L_Y^{23}(\lambda) L_X^{13}(\lambda-h^{(2)}) \\
& & \qquad \times\,  T_1 T_2 J^{12}(\lambda+h^{(2)}+h^{(3)})^{-1}
    x^{12}(\lambda+\hat h^{(1)}+\hat h^{(2)})\\
&=& \bar J^{12}x^{2}(\lambda-h^{(3)})^{-1}L_Y^{23}(\lambda)
    x^1(\lambda-h^{(2)}-h^{(3)})^{-1} L_X^{13}(\lambda-h^{(2)}) \\
& & \qquad \times\,  T_1 T_2 x^2(\lambda+\hat{h}^{(1)}+\hat{h}^{(2)})
    x^1(\lambda+\hat{h}^{(1)})(\bar J^{12})^{-1}\\
&=& \bar J^{12}x^{2}(\lambda-h^{(3)})^{-1}L_Y^{23}(\lambda){T_2}
    x^2(\lambda+\hat{h}^{(2)})x^1(\lambda-h^{(3)})^{-1} L_X^{13}(\lambda)\\
& & \qquad \times\,  {T_1} x^1(\lambda+\hat{h}^{(1)})(\bar J^{12})^{-1}\\
&=& \bar J^{12}\,\bar L_Y^{23}\, \bar L_X^{13}\,(\bar J^{12})^{-1},
\end{eqnarray*}
as desired.
\end{proof}

The functor ${\mathcal F}$ is quite nontrivial. For example, 
let $V$ be the trivial representation $(V=\C,L=1)$. Then 
${\mathcal F}(V)=F(A^*)$, with the action 
of $(H^J)^*$ given by:
$$
L_X=x^1(\lambda)^{-1}{T_1}x^1(\lambda+\hat{h}^{(1)}).
$$
This is typically a very nontrivial representation 
of $(H^J)^*$ by scalar difference operators on $A^*$. 
For example, for the vertex-IRF transformation in statistical mechanics
this very interesting
representation was considered by Hasegawa; it is connected to the 
theory of integrable systems, see \cite{KZ}. 

\end{section}


\begin{section}
{On classification and examples of   vertex-IRF transformations}

Now let us consider the special case when  $\bar A=\{0\}$. 
We will be interested in finding   vertex-IRF transformations of the
twist $1\otimes 1$. That is, we want to find functions $x(\lambda)$
such that 
\begin{equation}
\label{one-one-x}
(1\otimes 1)^x = \Delta(x)(\lambda)\,x^2(\lambda)^{-1}\,
x^1(\lambda-h^{(2)})^{-1}
\end{equation} 
is a dynamical twist (i.e., is of zero weight). 

We note, first of all, that the set ${\mathcal S}$ of such functions 
carries an action of the group of gauge transformations ${\mathcal G}$
of dynamical twists $A^*\to H\otimes H$, since one can always compose
a  vertex-IRF transformation with a gauge transformation
and obtain a new   vertex-IRF transformation.
What we are really interested in is the description 
of the set of orbits, ${\mathcal S}/{\mathcal G}$. 

Below we formulate a theorem which describes all elements of 
${\mathcal S}$ which normalize ${\C}[A]$.
Let $N_A$ be the normalizer of ${\C}[A]$ in
the multiplicative group $H^\times$. We have a natural 
homomorphism from $N_A$
to the group of permutations of $A^*$,
$\pi: N_A\to \Perm(A^*)$, defined by the action of 
$N_A$ on the primitive idempotents of $\C[A]$. 

\begin{theorem} 
\label{when x is vertex IRF}
A function $x:A^*\to N_A$ is a   vertex-IRF transformation 
of $1\otimes 1$ if and only if 
\begin{equation}
\pi(x(\lambda))^{-1}(\mu) = f(\lambda)-f(\lambda-\mu),
\end{equation}
for a suitable bijective function $f: A^*\to A^*$. 
Furthermore, two such  vertex-IRF transformations $x_i$, $i=1,2$, 
are gauge equivalent if and only if they define the same permutations
$\pi(x_i(\lambda))$.
\end{theorem}
\begin{proof} 
Let $\pi(x(\lambda))^{-1}(\mu)=F_\lambda(\mu)$. 
It is easy to see that the condition for 
$x$ to be a vertex-IRF transformation is
$$
F_\lambda(a+b)=F_\lambda(b)+F_{\lambda-b}(a),
$$
for all $a,b\in A^*$.
Setting $b=0$ one gets  $F_\lambda(0)=0$. Also, setting $b=\lambda$, we find 
$$
F_\lambda(a+\lambda)=f(\lambda)-g(a),
$$
where $f(\lambda) = F_\lambda(\lambda)$ and $g(a) =F_0(a)$.
Putting in the last equation $a=-\lambda$, we find $g(a)=f(-a)$, i.e. 
$$
F_\lambda(a+\lambda)=f(\lambda)-f(-a), 
$$
as desired. The converse statement is straightforward.
\end{proof}

In the situation of this theorem, we will say that $x$ {\em realizes} 
the function $f$. 

\begin{remark} 
Theorem~\ref{when x is vertex IRF} shows that the 
set of vertex-IRF transformations normalizing $\C[A]$,
modulo gauge equivalences, is finite. 
They, however, do not exhaust all the variety of vertex-IRF transformations; 
in particular, the original vertex-IRF transformations coming from 
physics are not of this type. 
\end{remark}

For general $H$, it is not easy to find all functions 
$f : A^*\to A^*$ realized by vertex-IRF transformations.
Below we use the classification results of \cite{EN2} to derive
an answer in the case when $H=\C[G]$, where $G$ is a  finite group.

Recall \cite[Theorem 6.6]{EN2} that the gauge equivalence classes of dynamical
twists $J(\lambda) : A^* \to \C[G] \otimes \C[G]$ are in bijective
correspondence with isomorphism classes of {\em dynamical data} for $(G, A)$,
i.e., collections $(K,\, \{V_\lambda \}_{\lambda\in A^*})$, where
$K$ is a subgroup of $G$ and $V_\lambda$ are irreducible projective
representations of $K$ such that $V_\lambda \otimes V_\mu^*$
is linear and
\begin{equation*}
\Ind_K^G(V_\lambda \otimes V_\mu^*) \cong \Ind_A^G(\lambda-\mu),
\end{equation*}
for all $\lambda,\mu\in A^*$ (see \cite[Section 4]{EN2} for a detailed
discussion of this notion).

\begin{proposition} 
\label{the one Pasha left me to prove}
The functions $f$ which are realized by vertex-IRF transformations
are those  for which 
\begin{equation}
\label{ind =ind}
\Ind_A^G(\lambda-\mu)  \cong \Ind_A^G(f(\lambda)-f(\mu))
\end{equation} 
for all $\lambda,\mu\in A^*$.
In particular, a dynamical twist is obtained from $1\otimes 1$ by 
a vertex-IRF transformation if and only if in the dynamical 
data corresponding to it, one has $K=A$. 
\end{proposition}
\begin{proof} 
If $f$ is realized by a vertex-IRF transformation $x$,
then for all $\lambda,\mu\in A^*$ the characters $\lambda-\mu$
and $f(\lambda)-f(\mu)$ are conjugate via $\Ad_{\pi(x(\lambda))^{-1}}$,
and therefore induce equivalent representations of $G$.

Conversely, if (\ref{ind =ind}) holds for a bijective function 
$f:A^* \to A^*$ then the collection 
$(A,\, \{\C_{f(\lambda)} \}_{\lambda\in A^*})$
is a dynamical datum for $(G, A)$. For every element $y\in \C[G]$
of weight $\mu$ the homomorphism 
$\Psi(\lambda, y) : \C_{f(\lambda)} \to \C_{f(\lambda-\mu)} \otimes \C[G]$
arising in the exchange construction \cite[Section 6]{EN2}, is given by
$1\mapsto 1\otimes x(\lambda)y$, where the element $x(\lambda)$ belongs
to the normalizer of $\C[A]$ in $\C[G]^\times$. In every $G$-module  
$x(\lambda)$ maps
elements of weight $\mu$ to elements of weight $f(\lambda)-f(\lambda-\mu)$,
for all $\lambda,\mu\in A^*$. Clearly, $x(\lambda)$ realizes $f$ and
the twist defined by the above action is equal to $(1\otimes 1)^x$.
\end{proof}

As an example, consider a class of vertex-IRF transformations in $\C[G]$ 
defined by elements $x$ which may be called ``quasi-grouplike elements''
(this example was considered in \cite{EN2}).

\begin{proposition}
\label{quasigr}
Let $g:A^*\times A^*\to G$ be a function, such that 
\begin{equation*}
(\lambda-\mu)\circ \Ad_{g(\lambda,\mu)} =f(\lambda)-f(\mu).
\end{equation*}
Then 
\begin{enumerate}
\item[(i)] 
$x(\lambda):=g(\lambda,\lambda-h^{(1)})$ is a vertex-IRF transformation 
of $1\otimes 1$ realizing $f$.
\item[(ii)] 
The transformed dynamical twist has the form
\begin{multline}
\label{J in terms of g}
J(\lambda)=:g(\lambda,\lambda-h^{(1)}-h^{(2)})g^{-1}(\lambda-h^{(2)},
\lambda-h^{(1)}-h^{(2)})\otimes \\
\otimes g(\lambda,\lambda-h^{(1)}-h^{(2)})g^{-1}(\lambda, \lambda-h^{(2)}):,
\end{multline}
\end{enumerate}
where the colons on both sides mean ``normal ordering'' : 
elements $h^{(1)}$ and $h^{(2)}$ are put to the extreme right, i.e.,
they are replaced by $\mu,\nu$ respectively when the twist acts on a vector 
$v\otimes w$ such that $v,w$ have weights $\mu,\nu$. 
\end{proposition}
\begin{proof}
The proof is contained in \cite[Example 6.10]{EN2}.
\end{proof}

The theory developed in this Section can be extended to the case when the 
finite group $A$ is replaced by an abelian Lie algebra ${\mathfrak a}$, 
along the lines  described in Section~2. 
We will not do it completely, but will just show 
how the generalization of the Proposition \ref{quasigr} allows one to quantize 
some classical dynamical $r$-matrices. 

Let $\g$ be a finite-dimensional Lie algebra over $\C$ with an
Abelian Lie subalgebra ${\mathfrak a}$. 
Let $G$ be the corresponding Lie group.
Let $f:{\mathfrak a}^*\to {\mathfrak a}^*$ be a meromorphic mapping.

We have the following analogue of Proposition \ref{quasigr}.
Let $g:{\mathfrak a}^*\times {\mathfrak a}^*\to G$ be a meromorphic function 
well defined at generic points of the diagonal
(in fact, it is only needed that $g$ be defined in the formal neighborhood 
of the diagonal).

\begin{proposition} 
\label{x is vertexIRF}
Assume that for $\lambda, \mu\in {\mathfrak a}^*$ we have
\begin{equation}
(\lambda-\mu)\circ  \Ad_{g(\lambda,\mu)} =f(\lambda)-f(\mu).
\end{equation}
Then 
\begin{enumerate}
\item[(i)] 
$x(\lambda):=g(\lambda,\lambda-\hbar 
h^{(1)})$ is a vertex-IRF transformation 
of $1\otimes 1$ realizing $f$,
i.e., $\pi(x(\lambda))^{-1}(\mu)=\hbar^{-1}(f(\lambda)-f(\lambda-\hbar\mu))$.
\item[(ii)] 
The transformed dynamical twist has the form
\begin{multline}
J(\lambda)=:g(\lambda,\lambda-\hbar(h^{(1)}+h^{(2)}))
g^{-1}(\lambda-\hbar  h^{(2)}, \lambda-\hbar(h^{(1)}+h^{(2)}))\otimes \\
\otimes  g(\lambda,\lambda-\hbar (h^{(1)}+h^{(2)}))g^{-1}(\lambda,
\lambda-\hbar h^{(2)}):.
\end{multline}
\end{enumerate}
\end{proposition}

Let us now calculate the quasi-classical limit of 
$J(\lambda)$. This means, having $J(\lambda)=1+\hbar \rho(\lambda)
+O(\hbar^2)$ to calculate the classical dynamical $r$-matrix
$r(\lambda):=\rho^{21}(\lambda)-\rho(\lambda)$. 
This $r$-matrix is given by the following proposition. 

\begin{proposition} 
Let $\{y_i\}$ be a basis of ${\mathfrak a}$ and $\{y^i\}$ the dual basis 
of ${\mathfrak a}^*$. Let $g(\lambda,\lambda)=\gamma(\lambda)$. 
Then 
\begin{equation}
r(\lambda)=\sum_i \, \frac{\partial \gamma(\lambda)}{\partial y^i}
\gamma^{-1}(\lambda) \wedge y_i
\end{equation}
\end{proposition}
\begin{proof}
The proof is by a direct calculation. 
\end{proof}

\begin{example}
\label{GL}
Consider $G=GL(n)\ltimes \C^n$, $\g=\text{Lie}(G)$, and 
${\mathfrak a}=\C^n$. Then for any function $f$ whose Jacobian 
does not vanish identically, we can define the adjoint 
$g(\lambda,\mu)^* : {\mathfrak a}^* \to {\mathfrak a}^*$ 
of $ g(\lambda,\mu)$ (which is regarded as an element of $GL(n) \subset G$)
by 
\begin{equation}
\label{g via f}
g(\lambda,\mu)^* = \sum_{m\ge 1}\,f^{(m)}(\lambda)(\mu-\lambda)^{m-1}/m!,
\end{equation}
which satisfies the required conditions (here $f^{(m)}(\lambda) :
({\mathfrak a}^*)^{\otimes m} \to {\mathfrak a}^*$ is a symmetric
linear map). In particular, 
$g(\lambda,\lambda)=f'(\lambda)^*  : {\mathfrak a}^* \to {\mathfrak a}^*$.
Thus, we have:
\begin{equation}
r(\lambda) = r_{n,f}(\lambda) = \sum_i\, 
 \frac{\partial f'(\lambda)^*}{\partial y^i}
(f'(\lambda)^*)^{-1} \wedge y_i.
\end{equation}
In particular, for $n=1$ we have a basis $\{X,Y\}$  of $\g$,
where $X$ generates $gl(1)$ and $Y$ generates ${\mathfrak a}=\C$, 
with $[XY]=Y$. Then 
\begin{equation}
r_{1,f}(\lambda)=\frac{f''(\lambda)}{f'(\lambda)} X\wedge Y.
\end{equation}
This classical $r$-matrix was considered by Xu in \cite{X}.

In this case, we have $g(\lambda,\mu)=
\left(\frac{f(\lambda)-f(\mu)}{\lambda-\mu}\right)^X$, and 
the twist has the form 
\begin{multline}
J(\lambda)= \\
= :\biggl( 
\frac{f(\lambda)-f(\lambda-\hbar(Y^{(1)}+Y^{(2)}))}{Y^{(1)}+Y^{(2)}} 
\frac{Y^{(1)}}{f(\lambda-\hbar Y^{(2)}) - f(\lambda-\hbar(Y^{(1)}+Y^{(2)}))}
\biggr)^{X^{(1)}}\times  \\
\times \biggl( \frac{f(\lambda)-f(\lambda-\hbar(Y^{(1)}+Y^{(2)}))}
{Y^{(1)}+Y^{(2)}}\frac{Y^{(2)}}{f(\lambda)-
f(\lambda-\hbar Y^{(2)})}
\biggr)^{X^{(2)}}:,
\end{multline}
where for any monomial $F$ in $X^{i}$ and $Y^{i}$, 
the expression $:F:$ means the monomial $F$ in which all the $Y$-factors have 
been moved to the right from the $X$-factors.
\end{example}

\begin{remark}
We note that we have given here a quantization 
for the dynamical $r$-matrices $r_{n,f}(\lambda)$. 
As far as we know,
such a quantization was previously unknown, 
even for $n=1$. To be more specific, the paper
\cite{X} constructs a quantization of a vast 
collection of skew-symmetric dynamical $r$-matrices, 
but does not cover the case of $r_{1,f}$,
except for some specific $f$. In fact, 
the above constructions were motivated by P.~Xu's 
question ``how to quantize $r_{1,f}$?''
\end{remark}

\begin{remark}
As a special case of $r_{1,f}$, one can consider $f(\lambda)=e^{\lambda}$. 
In this case, $r_{1,f}$ is a constant $r$-matrix $X\wedge Y$. 
This $r$-matrix can be quantized by the well known Jordanian twist \cite{GGS}, 
but here we have quantized it as a dynamical $r$-matrix, i.e., preserving the 
weight zero condition. It is easy to see that 
the quantization we constructed is also constant (i.e., is a 
usual twist), but it has zero weight with respect to $Y$ and therefore 
is more complicated than the usual Jordanian twist. 
The existence of such ``zero weight'' quantization follows from \cite{EK},
but here we wrote it explicitly.
\end{remark}

\end{section}

\begin{section}
{Completely degenerate dynamical $r$-matrices}

Let $\g$ be a finite-dimensional Lie algebra, and $\mathfrak{a}
\subset \g$ be an Abelian Lie subalgebra.
In the following, by functions on $\a^*$ we will mean 
holomorphic functions defined near $0\in \a^*$. 

Let $\{ y_i\}_{i=1}^n$ be a basis of $\a$ and 
$\{ y^i\}_{i=1}^n$ be the dual basis of $\a^*$.
Let $r : \a^* \to \wedge^2 \g$ be a skew-symmetric
classical dynamical $r$-matrix, 
i.e., a zero weight solution of the classical
dynamical Yang-Baxter equation

\begin{equation}
\label{E:3}
\begin{split}
\sum_i &\left(y_i^{(1)} 
\frac{\partial r^{23}(\lambda)}{\partial y^i}-y_i^{(2)} 
\frac{\partial r^{13}(\lambda)}{\partial y^i} + y_i^{(3)} 
\frac{\partial r^{12}(\lambda)}{\partial y^i}\right) +\\
&[r^{12}(\lambda),r^{13}(\lambda)]+[r^{12}(\lambda),r^{23}(\lambda)]+
[r^{13}(\lambda),r^{23}(\lambda)]=0.
\end{split}
\end{equation}
 
For basic facts about such $r$-matrices
and their quantization, see survey \cite{ES2}.

Recall \cite{X} that $r$ is said to be non-degenerate if 
the projection $\tilde r(\lambda)$ of the element $r(\lambda)$ to 
$\wedge^2 (\g / \a)$ is non-degenerate for some $\lambda$. 
It is shown in \cite{X} that such an $r$-matrix can be quantized.

This motivates the following definition. 

\begin{definition}
We will say that $r$ is {\em completely
degenerate} if the induced map $\tilde{r} : \a^* \to \wedge^2 
(\g / \a)$ is zero. 
\end{definition}

In this Section, we classify all completely degenerate dynamical $r$-matrices, 
and use ideas of Section 3 to show that any such
$r$-matrix can be quantized. This generalizes Example~\ref{GL}.

\begin{remark} 
We note that if $\a$ has an invariant complement in $\g$ then 
it follows from \cite{ES3} that any completely degenerate dynamical 
$r$-matrix $r : \a^* \to \wedge^2 \g$
is gauge equivalent to zero. However, if such a complement 
does not exist then the theory of \cite{ES3} does not apply, 
and in particular there are many nontrivial examples of 
completely degenerate $r$-matrices (for instance, $r_{n,f}$ 
considered in the previous section).
\end{remark}

Before stating the result, we need to introduce some notation. 
Let us denote by $\n(\a) = \{ x\in \g \mid [x,\a] \subset
\a \}$ the normalizer of $\a$ in $\g$, and let
 $N(\a)$ be the simply connected Lie group 
corresponding to $\n(\a)$. We have a natural homomorphism 
$N(\a)\to GL(\a)$, which we will denote by $\gamma\to\bar\gamma$. 

\begin{theorem} 
\label{main theorem}
\begin{enumerate}
\item[(i)] 
Let $\gamma: \a^*\to N(\a)$ be a  function, such that $\bar\gamma(\lambda)^*
=f'(\lambda)$  for some function $f:\a^*\to \a^*$. Then the function 
\begin{equation}
\label{completely degenerate : formula}
r_\gamma(\lambda) =  \sum_i\, 
\frac{\partial {\gamma}(\lambda)}{\partial y^i} {\gamma}(\lambda)^{-1}
\wedge y_i, 
\end{equation}
is a completely degenerate dynamical $r$-matrix.
\item[(ii)] 
Every completely degenerate skew-symmetric classical dynamical $r$-matrix
is equal to (\ref{completely degenerate : formula})
for some function $\gamma : \a^* \to N(\a)$.
Moreover, $r_{\gamma_1}$ is gauge equivalent to $r_{\gamma_2}$ 
in the sense of \cite[p.3]{ES3} if and only if $\bar\gamma_1=\bar\gamma_2$.
\item[(iii)] Every completely degenerate   
skew-symmetric classical dynamical $r$-matrix
can be quantized.
\end{enumerate}
\end{theorem}

The rest of the Section is the proof of the theorem.

\begin{proof}
Let us first prove (i). 
Let $r$ be given by (\ref{completely degenerate : formula}).
First of all, it is easy to check that 
$r$ has zero weight. Namely, the zero weight condition
coincides with the cross-derivative 
condition for $\bar\gamma(\lambda)^*$, which is equivalent to 
the condition that $\bar\gamma(\lambda)^*=f'(\lambda)$ for some $f$. 
It is also clear that $\tilde r=0$. 
So it remains to check the classical Yang-Baxter equation for $r$. 

Let us write 
\begin{equation}\label{r=py}
r(\lambda) = \sum_i \, p_i(\lambda) \wedge y_i
\end{equation}
with $p_i(\lambda)=
\frac{\partial {\gamma}(\lambda)}{\partial y^i} {\gamma}(\lambda)^{-1}
\in \g$. 
Then it is straightforward to check that the zero weight condition for
$r$ translates into
\begin{equation}\label{zw}
\sum_i\, y_i\otimes [y,p_i(\lambda)] = \sum_i\, [y,p_i(\lambda)] 
\otimes  y_i,
\end{equation}
for all $y\in \a$. Using this identity we can write the 
classical dynamical Yang-Baxter equation for $r$ in terms of functions
$p_i(\lambda)$ as
\begin{multline}
\label{cde for y and p}
\sum_{ij} \, \left(\frac{\partial p_j(\lambda)}{\partial y^i}
- \frac{\partial p_i(\lambda)}{\partial y^j} +[p_i(\lambda), p_j(\lambda)] 
\right) \otimes y_i \otimes y_j \\
 - y_i \otimes \left (\frac{\partial p_j(\lambda)}{\partial y^i}
- \frac{\partial p_i(\lambda)}{\partial y^j} +[p_i(\lambda), p_j(\lambda)]
\right) \otimes y_j \\
+ y_i \otimes y_j  \otimes  
\left (\frac{\partial p_j(\lambda)}{\partial y^i}
- \frac{\partial p_i(\lambda)}{\partial y^j} +[p_i(\lambda), p_j(\lambda)]
\right)=0.
\end{multline}
This equation is clearly satisfied, since 
$p_i=\frac{\partial {\gamma}(\lambda)}{\partial y^i} {\gamma}(\lambda)^{-1}$.
So statement (i) is proved.

Let us now prove (ii). Let $r$ be any completely degenerate 
dynamical $r$-matrix. 

\begin{lemma}
\label{where r lives}
We have $r(\lambda) \in \n(\a) \wedge \a$.
\end{lemma}
\begin{proof}
Clearly, $r(\lambda)\in \g \wedge \a$. 
Write $r$ in the form (\ref{r=py}).
Since $r$ has zero weight, 
 the projection
\begin{equation}
\hat{r}(\lambda) = \sum_i \, \hat{p}_i(\lambda) \otimes y_i,
\quad \text { where } \quad \hat p_i(\lambda)\in \g/\a,
\end{equation}
of $r(\lambda)$ on $(\g/\a)\otimes \a$ also has zero weight. 
This implies that $[y_j, \hat{p}_i(\lambda)] =0$ in $\g/\a$
for all $i,j$, therefore $\hat{p}_i(\lambda) \in \n(\a)/\a$,
and hence $p_i(\lambda) \in \n(\a)$.
\end{proof}

Now, as we explained above, 
the classical dynamical Yang-Baxter equation 
for $r$ and the zero weight condition
reduce to equations (\ref{cde for y and p})
and (\ref{zw}) respectively. In particular, we have 
\begin{equation}
\label{p-equation}
\frac{\partial \hat{p}_j(\lambda)}{\partial y^i}
- \frac{\partial \hat{p}_i(\lambda)}{\partial y^j} +
[\hat{p}_i(\lambda), \hat{p}_j(\lambda)] = 0.
\end{equation}
This implies that the differential equations
\begin{equation}
\frac{\partial \hat{\gamma}(\lambda)}{\partial y^i} = \hat{p}_i(\lambda)
\hat{\gamma}(\lambda),
\end{equation}
are compatible and have a unique solution $\hat\gamma: \a^*\to 
N(\a)/\exp(\a)$ with the initial condition $\hat{\gamma}(0)=1$. 

Let $\gamma : \a^* \to N(\a)$ be an arbitrary lift of  $\hat{\gamma}$. Then 
\begin{equation}
\label{r+gauge}
r(\lambda) =  \sum_i\, 
\frac{\partial {\gamma}(\lambda)}{\partial y^i} {\gamma}(\lambda)^{-1}
\wedge y_i + \sum_{ij}\, C_{ij}(\lambda) y_i \wedge y_j, 
\end{equation}
where the second sum is a  $2$-form.  It follows from
the classical dynamical Yang-Baxter equation that $C$ is a closed form.
Therefore, $r$ is gauge equivalent to 
(\ref{completely degenerate : formula}), as desired. 
This proves the first statement of (ii). 

To prove the second statement of (ii),
it is sufficient to recall that gauge transformations 
are functions from $\a^*$ to the centralizer $Z(\a)$ of $\a$ in $N(\a)$.
Now, the action of such a function $g$ on $r_\gamma$ is given by 
$\gamma\to g\gamma$. This implies easily the second statement of (ii),
since the image of $N(\a)$ in $GL(\a)$ is exactly $N(\a)/Z(\a)$. 

Finally, let us prove (iii).  
First of all, as we already mentioned in the proof of (i), 
the zero weight property of $r$ implies
that $\bar\gamma(\lambda)^* = f'(\lambda)$ for some function $f:
\a^* \to\a^*$. Now, for $\lambda\in \a^*$ define the element
$$
x(\lambda)=\sum_{m\ge 1}\frac{(-1)^{m-1}}{m!}
\gamma^{(m-1)}(\lambda)(\hbar h^{(1)})^{m-1}
$$
of the  group $\tilde H^\times$ for 
$H=U(\n(\a))[[\hbar]]$. Then it is easy to check that 
$x(\lambda)h^{(1)}=
\hbar^{-1}(f(\lambda)-f(\lambda-\hbar h^{(1)})x(\lambda)$.
Therefore, 
by Section 3, $x(\lambda)$ is a vertex-IRF transformation of $1\otimes 1$.
Hence $J(\lambda)=\Delta(x(\lambda))x^2(\lambda)^{-1}x^1(\lambda-h^{(2)})^{-1}$
is a dynamical twist. It is easy to check directly that 
the quasi-classical limit of $J(\lambda)$ is $r(\lambda)$. 
The theorem is proved.
\end{proof}

\begin{remark}
We warn the reader that $x(\lambda)$ is not quasi-grouplike, 
and hence for the examples of Section 3, 
the quantization constructed here coincides with the quantization by 
``quasi-grouplike elements'' of Section 3 only up to gauge 
transformations. 
\end{remark}

\end{section}


\bibliographystyle{ams-alpha}

\end{document}